# Effectiveness in RPL, with Applications to Continuous Logic

Farzad Didehvar *
didehvar@aut.ac.ir

Kaveh Ghasemloo [†]
[first name]@cs.toronto.edu

Massoud Pourmahdian[*‡§]
pourmahd@ipm.ir

June 11, 2009

## Abstract

In this paper, we introduce a foundation for computable model theory of rational Pavelka logic (an extension of Łukasiewicz logic) and continuous logic, and prove effective versions of some related theorems in model theory. We show how to reduce continuous logic to rational Pavelka logic. We also define notions of computability and decidability of a model for logics with computable, but uncountable, set of truth values; we show that provability degree of a formula with respect to a linear theory is computable, and use this to carry out an effective Henkin construction. Therefore, for any effectively given consistent linear theory in continuous logic, we effectively produce its decidable model. This is the best possible, since we show that the computable model theory of continuous logic is an extension of computable model theory of classical logic. We conclude with noting that the unique separable model of a separably categorical and computably axiomatizable theory (such as that of a probability space or an $L^p$ Banach lattice) is decidable.

---
*Department of Mathematics and Computer Science, Amirkabir University of Technology, Tehran, Iran
[†]Department of Computer Sciences, University of Toronto, Toronto, ON, Canada
[‡]School of Mathematics, Institute for Research in Fundamental Sciences (IPM), Tehran, Iran
[§]Partially supported by a grant from IPM.

# 1 Introduction

The school of computable mathematics, which is one of the major constructive approaches toward the subject, studies mathematical structures from the computability theoretic viewpoint. The importance of this school increased significantly after the realization of digital computers, and after many well established subjects in mathematics were restudied through the theme of computability; classical model theory (model theory of classical logic) was not an exception. [Har98, EGNR98, Ers98a, Ers98b, Mil99] revisit many classical concepts and theorems in classical model theory under the light of computability. Extending this approach to various non-classical logics has been carried out, e.g., in [GN04].

First order logic serves well for a model theoretic study of structures in combinatorics and algebra. On the other hand, it is not very fruitful to use classical first order logic for model theoretic study of structures in analysis. Therefore, a direct approach to use computability over uncountable domains for generalizing results is not viable. Using a non-classical logic can be more fruitful for model theoretic study of analytic structures. *Continuous logic* (**CL**), introduced by Chang and Keisler in [CK66], seems to be a natural choice for this endeavor. Working with **CL** has become easier to a considerable extent due to [BU] and [BBHU].

In this paper, we study *computable continuous model theory* (CCMT), i.e. *model theory of continuous logic computably*.

Use of a non-classical logic, such as modal or continuous logic, generates considerable complexity in effectivizing arguments, even in the most basic model theoretic constructions; one such example is in establishing the completeness theorem via a Henkin construction as seen through comparing the case of classical logic [Mil99] with the case of modal one [GN04]. Not every formulation of classical completeness theorem is correct for continuous logic; furthermore, similar to modal logics, neither the deduction theorem nor the principle of excluded middle (PEM) holds for continuous logic.

Another main obstacle in extending computable model theory to continuous logic is the uncountability of truth values involved. The concepts in classical computability theory were mainly developed for countable domains such as natural numbers, and structures that can be represented by them. These concepts are adopted for the study of effectiveness in countable structures, not uncountable ones. Again, combinatorics and algebra fit well in this framework, where many of interesting structures occur naturally in the finite or countable cases. On the other hand, we have analysis, wherein natural interesting structures are essentially based on real numbers and thus contain uncountable concepts, and which does not embed well in the classical computability theory. There has been some interest in computability on uncountable domains since late 1950s but this topic remained mostly out of the mainstream of research in computability theory till recently. The need for robust and reliable computer software dealing with computation over real numbers, or in analytic or geometric structures, has changed the situation. Currently there are many proposed foundations, which are, to various extents, incompatible with each other.[1] The main idea behind many of them is the possibility of approximating a result with arbitrary precision given good enough approximation(s) for the necessary inputs. We use *continuous domain* theoretic approach to present our results, as it makes the presentation simpler, although translating them to TTE or Scott-Ershov $\omega$-algebraic approach is straightforward.

Here, we initiate the study of computable model theory for continuous logic. We believe that CCMT fills the place of the question mark in the following analogy, thereby offering a satisfactory and unifying framework to the model theoretic, effective study of analytic structures:

$$\frac{\text{Computable (Classical) Model Theory}}{\text{Computable Algebra}} \approx \frac{?}{\text{Computable Analysis}}$$

The first step in this direction is defining primary concepts and adapting basic theorems such as Henkin's model construction technique to our framework.

## 2 Preliminaries

### 2.1 Łukasiewicz Logic, Rational Pavelka Logic

We use **LL∀** to denote *first order Łukasiewicz logic*. Łukasiewicz logic is the only many-valued logic given by a continuous t-norm and a continuous residuum. **RPL∀** is *first order rational Pavelka logic*, the extension of **LL∀** by propositional truth values for rational numbers. Following [MOG]'s notions, the logical operations are: $\oplus, \odot, \wedge, \vee, \rightarrow, \leftrightarrow, \neg, \forall, \exists, \bot$. The unit interval $[0,1]$ is our standard lattice for the evaluation of operators; $\oplus$ is interpreted by the continuous t-norm $x \dotplus y := \min\{x+y, 1\}$; $\odot$ by $x * y := \max\{1-x-y, 0\}$; $\wedge$ and $\vee$ are interpreted by max and min over $[0,1]$, respectively; $x \rightarrow y$ by the $y \dotminus x := \max\{y-x, 0\}$; $\neg$ by $x \mapsto (1-x)$; $\leftrightarrow$ by $(x, y) \mapsto |x-y|$; $\forall$ by sup; and $\exists$ by inf. We will use $\bigsqcup$ for sup, interpreting both $\wedge$ and $\forall$; and $\bigsqcap$ for inf, interpreting both $\vee$ and $\exists$. The truth value of $\bot$ is 1. A propositional constant $\bar{r}$ has truth value $r$, for $r$ in $\mathbb{Q}_{[0,1]}$ (the set of rationals in $[0,1]$).[2] The following conventions are used: $x, y, z, \ldots$ represent variables; $a, b, c, \ldots$ represent constants; $f, g, h, \ldots$ represent functions; $r, s, t, \ldots$ represent terms; $P, Q, R, \ldots$ represent relations; and $\varphi, \psi, \sigma, \ldots$ represent formulas. $n\varphi$ and $\varphi^n$ abbreviate $\underbrace{\varphi \oplus \ldots \oplus \varphi}_{\text{n-times}}$ and $\underbrace{\varphi \odot \ldots \odot \varphi}_{\text{n-times}}$.

---

[1] E.g., TTE: [Wei00], Banach/Mazur, Grzegorczyk [Grz55, Grz57], Pour-El/Richard [PER89, PE99], Ko [Ko91, Ko98]; Continuous Domains: [EH98, ES99b, ES99a, Eda97]; Scott-Ershov Algebraic Domains: [SHT07, Bla00, Bla97b, Bla97a]; Markov School, Real-RAM: [BSS89, BCSS96, BCSS98], Sharp Filters [KW98, KW99]. For a comparison of these, see [Wei00, §9].

[2] Notice that we consider 0 as the highest degree of truth (the least degree of falsity), and 1 as the least. Although this is not the usual practice, see [Haj98] for example, it simplifies the interpretation of operators. This will be helpful when we consider *continuous logic* later. Also, the isomorphism $x \mapsto (1-x)$ can be used to get [Haj98]'s notions.

A *standard structure* for a *language* $L$ is defined similar to the classical bi-valued Boolean logic, but relations take values in $[0, 1]$. An *$M$-evaluation* $v$ is a mapping from set of variables to the universe of $M$. The *truth degree* of $\varphi$ in $(M, v)$, $||\varphi||_{M,v}$, is defined as the truth value of its interpretation, defined inductively from atomic formulas using interpretation of logical operators. $||\varphi||_M$ and $||\varphi||$ abbreviate, respectively, $\sqcup\{||\varphi||_{M,v} : v \text{ an } M-evaluation\}$ and $\sqcup\{||\varphi||_M : M \text{ a standard structure}\}$. $\varphi$ is a *tautology* iff $||\varphi|| = 0$. A *model* of a *theory* $T$ in the language $L$, is a standard structure for $L$, $M$, and an evaluation $v$, so that the truth degree of any formula of $T$ in $(M, v)$ is 0; we write $(M, v) \vDash T$ in this case. If for every $v$, $(M, v) \vDash T$, we write $M \vDash T$. $||\varphi||_T$ is defined as $\sqcup\{||\varphi||_M : M \vDash T\}$.

We take $\rightarrow$, $\bot$, and $\forall$ as primitive connectives and define the rest as follows:

$\neg$: $\neg\varphi := \varphi \rightarrow \bot$,

$\oplus$: $\varphi \oplus \psi := \neg(\varphi \rightarrow \neg\psi)$,

$\odot$: $\varphi \odot \psi := \neg\varphi \rightarrow \psi$,

$\wedge$: $\varphi \wedge \psi := \varphi \oplus (\varphi \rightarrow \psi)$,

$\vee$: $\varphi \vee \psi := (\varphi \rightarrow \psi) \rightarrow \psi$,

$\leftrightarrow$: $\varphi \leftrightarrow \psi := (\varphi \rightarrow \psi) \oplus (\psi \rightarrow \varphi)$,

$\exists$: $\exists x\, \varphi(x) := \neg\forall x\, \neg\varphi(x)$.

The axioms and rules of **LL∀** from [MOG, p. 94], [Haj98, 3.1.3]:

## LL∀

**L1:** $\varphi \rightarrow (\psi \rightarrow \varphi)$,

**L2:** $(\varphi \rightarrow \psi) \rightarrow ((\psi \rightarrow \sigma) \rightarrow (\varphi \rightarrow \sigma))$,

**L3:** $((\varphi \rightarrow \psi) \rightarrow \psi) \rightarrow ((\psi \rightarrow \varphi) \rightarrow \varphi)$,

**L4:** $((\varphi \rightarrow \bot) \rightarrow (\psi \rightarrow \bot)) \rightarrow (\psi \rightarrow \varphi)$,

**∀1:** $\forall x\, \varphi \rightarrow \varphi[t/x]$, where no free variable of $t$ becomes bounded by substitution for $x$ in $\varphi$,

**∀2:** $\forall x\, (\varphi \rightarrow \psi) \rightarrow (\varphi \rightarrow \forall x\, \psi)$, where $x$ is not free in $\varphi$,

**mp:** $\frac{\varphi \quad \varphi \rightarrow \psi}{\psi}$,

**gen:** $\frac{\varphi}{\forall x\, \varphi}$.

And **RPL∀** is given by adding:

## RPL∀

**R:** $(\bar{r} \to \bar{s}) \leftrightarrow \overline{s \mathbin{\dot{-}} r}$, for $r$ and $s$ in $\mathbb{Q}_{[0,1]}$.

**RPL∀** is a definable extension of **LL∀**. For $r \in \mathbb{Q}_{[0,1]}$, let $\varphi_r$ be new atomic propositions, and consider the following axioms:

$$\varphi_0 \leftrightarrow (\varphi_0 \odot \neg\varphi_0),\ \varphi_1 \leftrightarrow (\varphi_1 \oplus \neg\varphi_1),\ \varphi_{\frac{1}{2}} \leftrightarrow \neg\varphi_{\frac{1}{2}},\ n\varphi_{\frac{1}{2n}} \leftrightarrow \varphi_{\frac{1}{2}},\ \varphi_{\frac{m}{n}} \leftrightarrow 2m\varphi_{\frac{1}{2n}}.$$

The completeness theorem of **LL∀** holds when restricted to standard models, which are also models of **RPL∀**; therefore adding these axioms will not change properties of theories such as consistency.

Following [Haj98], we do not consider *equality* as a logical relation, and do *not* interpret it by *identity*. This means that there are distinct elements which are equal with the truth degree 0. We do not have any set of axioms forcing the equality, which is a congruence relation, to be the true identity. Theoretically, one can create a model that interprets equality as identity in a model, by taking a quotient over the crisp equality relation of the model, but as we will see, this process is not always effective.

Following is the list of axioms for fuzzy equality, $\approx$, which is called *similarity*:

### Similarity

$S1$(**Reflexivity**): $x \approx x$,

$S2$(**Symmetry**): $x \approx y \to y \approx x$,

$S3$(**Transitivity**): $(x \approx y \oplus y \approx z) \to x \approx z$.

It is easy to check that any interpretation of the similarity relation, satisfying above axioms, is actually a pseudo-metric. It is sometimes useful to require similarity to satisfy the following axioms:

### Congruence

$S4_R$: $\overrightarrow{x} \approx \overrightarrow{y} \to (R \leftrightarrow R[\overrightarrow{x}/\overrightarrow{y}])$, where $\overrightarrow{x} = (x_1, ..., x_n)$, $\overrightarrow{y} = (y_1, ..., y_n)$, and $n$ is the arity of $R$, $\overrightarrow{x} \approx \overrightarrow{y}$ abbreviating $x_1 \approx y_1 \oplus \cdots \oplus x_n \approx y_n$,

$S4_f$: $\overrightarrow{x} \approx \overrightarrow{y} \to (f(\overrightarrow{x}) \approx f(\overrightarrow{y}))$, with the same conditions as above.

In every structure that $S4$ holds, all interpretations will become 1-Lipchitz. Such structures are called *extensional* [3]. $S4$ is not assumed below, unless explicitly stated.

We write $T \vdash \varphi$ when there is a proof for $\varphi$ from $T$ in **RPL∀**. $\varphi$ is *valid* if $\vdash \varphi$. The *provability degree* of $\varphi$ over $T$, $|\varphi|_T$ is defined as $\sqcap\{r \in \mathbb{Q}_{[0,1]} : T \vdash \bar{r} \to \varphi\}$. We say $T$ is *consistent* if $T \nvdash \bot$. We say $T$ is *strongly consistent* if $\forall r \in \mathbb{Q}_{(0,1]} \, T \nvdash \bar{r}$. In a standard model $M$, and for all $r \in \mathbb{Q}_{(0,1]}$, we have $M \nvDash \bar{r}$. This is the motivation behind the definition. Every satisfiable theory is strongly consistent. This condition is necessary to prove the existence of a standard model for a theory, but it is easy to see this is equivalent to consistency (the less obvious direction holds by taking a natural $n$ such that $nr \geq 1$ and showing $T \vdash n\bar{r}$). $T$ is *complete* iff for any sentence $\varphi$, either $T \vdash \varphi$ or $T \vdash \neg\varphi$. $T$ is *linear complete* if for any two sentences $\varphi$ and $\psi$, either $T \vdash \varphi \to \psi$ or $T \vdash \psi \to \varphi$.[4] $T$ is *Henkin* if for any sentence of the form $\exists x\, \varphi$, if $T \vdash \exists x\, \varphi$, then there is a constant $c$ such that $T \vdash \varphi[c/x]$.

---
[3] See [Haj98, p.142]

For a consistent, linear complete $T$, provability degree is equal to $\sqcup \{r \in \mathbb{Q}_{[0,1]} : T \vdash \varphi \to \bar{r}\}$ [Haj98, Theorem 3.3.8 part 1]. Provability degree commutes with logical connectives [Haj98, Theorem 3.3.8 part 2]. It also commutes with quantifiers over the constants of language, if the theory is Henkin [Haj98, Theorem 5.2.6 part 2].

*Deduction Theorem* is not valid for **RPL∀**, but we have the following weak form [Haj98]:

**Lemma 2.1.** *If $T, \varphi \vdash \psi$, then for some natural $n$, $T \vdash n\varphi \to \psi$.*[5]

**RPL∀** is *sound*: if $T \vdash \varphi$, then $T \vDash \varphi$, with respect to models. **RPL∀** also satisfies completeness theorem: $||\varphi||_T = |\varphi|_T$.[6]

In **RPL∀**, $T \vDash \varphi$ and $T \vdash \varphi$ are not equivalent. The following lemma shows their relationship.

**Lemma 2.2.** *For all $q > 0$, $T \vdash \varphi \dotdiv q$ iff $T \vDash \varphi$.*

Hence, there are two options for axiomatizability. Axiomatizability using '$\vdash$' leads to axiomatizability using '$\vDash$'. Thus the latter will be called *weak axiomatizability*.

We list some theorems from [Haj98], which we will use.

**Lemma 2.3.** *If $T \nvdash \varphi \to \psi$, then $T \cup \{\psi \to \varphi\}$ is consistent.*

**Lemma 2.4.** *If every finite $T_0 \subseteq T$ has a model, then $T$ has a model.*

**Lemma 2.5.** *$T$ is linear complete iff for each pair $(\varphi, \psi)$ of sentences, if $T \vdash \varphi \vee \psi$, then either $T \vdash \varphi$ or $T \vdash \psi$.*

**Lemma 2.6.** *The logical connectives, $\to, \forall$ commutate with provability degree for consistent, linear complete, Henkin theories; $|\varphi_2|_T \dotdiv |\varphi_1|_T = |\varphi_1 \to \varphi_2|_T$, $\sqcup_{a \in C} |\psi[x/a]|_T = |\forall x\, \psi|_T$.*

**Lemma 2.7.** $\vdash (\varphi \to \psi)^n \vee (\psi \to \varphi)^n$.

---

[4][Haj98] calls this *complete*. *Linear complete* seems more appropriate, because the condition forces the *Lindenbaum Algebra* of theory to be a linear order, and it will not conflict with completeness in classical sense.

[5]See [Haj98, p.43].

[6] Let $T = \{np \to q : n \in \omega\} \cup \{\neg p \to q\}$. Then $T \vDash q$, but $T \nvdash q$! Therefore $T \vDash \varphi \Rightarrow T \vdash \varphi$ does not hold.

There are a few notions for completeness of a theory, to avoid confusion, we list their definitions.

1. Classical completeness: for all $\varphi$ either $T \vdash \varphi$ or $T \vdash \neg\varphi$,

2. Linear completeness: for all $\varphi$ and $\psi$, either $T \vdash \varphi \to \psi$ or $T \vdash \psi \to \varphi$,

3. Semantically linear completeness: for all $\varphi$ and $\psi$, either $T \vDash \varphi \to \psi$ or $T \vDash \psi \to \varphi$.

## 2.2 Continuous Logic

We follow [BU] and [BBHU]. *Continuous Logic* (**CL**) is a truth functional logic. The set of truth values is $[0, 1]$. We call a set of continuous functions over $[0, 1]$ a *full system of connectives*, if it contains projections, is closed under composition, and is dense in compact-open (i.e. uniform convergence) topology over the set of all continuous functions over $[0, 1]$. By Stone-Weierstrass Theorem [BU, §1.6], $\{\dotminus\} \cup \{\bar{q} : q \in \mathbb{Q}_{[0,1]}\}$ is a full system. The natural reinterpretation of classical quantifiers in this continuous setting are sup and inf [BU, p. 6]. We will use $\dotminus$, $\bar{q}$ for $q \in \mathbb{Q}_{[0,1]}$, and sup as our connectives. Other connectives can be defined from these as we did for **LL**∀ above.

We call a function $\delta : (0, 1] \to (0, 1]$ a *modulus of continuity*. Assume that $(X_1, d_1)$ and $(X_2, d_2)$ are two (pseudo-)metric spaces, and $f$ a function from $X_1$ to $X_2$. We say $f$ is *uniformly continuous with respect to $\delta$*, iff $\forall \epsilon > 0 \forall x, y \in X_1(d_1(x, y) < \delta(\epsilon) \to d_2(f(x), f(y)) \leq \epsilon)$. A *non-metric continuous signature* is a set of function and relation symbols with their arities. A *(metric) continuous signature* is a non-metric continuous signature with a distinguished binary relation $d$ that has, for each $s$ of the non-metric continuous signature and for all $i < n_s$, a uniform continuity modulus $\delta_{s,i}$.

Let $L$ be a continuous signature. A *continuous pre-structure* is a set $M$, a pseudo-metric $d^M$ over $M$, where for each function symbol $f \in L$, there is a function $f^M : M^{n_f} \to M$ that is uniformly continuous with respect to $\delta_{f,i}$ $(i < n_F)$; and where for each relation symbol $R \in L$, there is a function $R^M : M^{n_R} \to [0, 1]$ that is uniformly continuous with respect to $\delta_{R,i}$ $(i < n_R)$. A *continuous structure* is a pre-structure where $d$ is complete metric. From any pre-structure, we can build an elementary equivalent structure by first taking quotient (with respect to distance-zero points), and then completing it (with respect to its metric).

The axioms and rules of **CL** are:

$$\textbf{CL}$$

*C1:* $(\varphi \dotminus \psi) \dotminus \varphi$,

*C2:* $((\sigma \dotminus \varphi) \dotminus (\sigma \dotminus \psi)) \dotminus (\psi \dotminus \varphi)$,

*C3:* $(\varphi \dotminus (\varphi \dotminus \psi)) \dotminus (\psi \dotminus (\psi \dotminus \varphi))$,

$C4$: $(\varphi \mathbin{\dot{-}} \psi) \mathbin{\dot{-}} ((\bar{1} \mathbin{\dot{-}} \psi) \mathbin{\dot{-}} (\bar{1} \mathbin{\dot{-}} \varphi))$,

sup 1: $\varphi[t/x] \mathbin{\dot{-}} \sup_x \varphi$, such that no free variable of $t$ becomes bounded by substitution for $x$ in $\varphi$, sup 2: $(\sup_x \psi \mathbin{\dot{-}} \varphi) \mathbin{\dot{-}} \sup_x(\psi \mathbin{\dot{-}} \varphi)$, where $x$ is not free in $\varphi$,

$R1$: $(\bar{r} \mathbin{\dot{-}} \bar{s}) \mathbin{\dot{-}} \overline{r \mathbin{\dot{-}} s}$, for $r$ and $s$ in $\mathbb{Q}_{[0,1]}$,

$R2$: $\overline{r \mathbin{\dot{-}} s} \mathbin{\dot{-}} (\bar{r} \mathbin{\dot{-}} \bar{s})$, for $r$ and $s$ in $\mathbb{Q}_{[0,1]}$,

$SM1$: $d(x,x)$,

$SM2$: $d(x,y) \mathbin{\dot{-}} d(y,x)$,

$SM3$: $(d(x,z) \mathbin{\dot{-}} d(x,y)) \mathbin{\dot{-}} d(y,z)$,

$UL_R$: $(q \mathbin{\dot{-}} d(x,y)) \vee ((R(\bar{a},x,\bar{b}) \mathbin{\dot{-}} R(\bar{a},y,\bar{b})) \mathbin{\dot{-}} r)$, where $\epsilon, r, q \in \mathbb{Q}_{[0,1]}$, $r > \epsilon$, $q < \delta_{R,i}(\epsilon)$,

$UL_f$: $(q \mathbin{\dot{-}} d(x,y)) \vee (d(f(\bar{a},x,\bar{b}), f(\bar{a},y,\bar{b})) \mathbin{\dot{-}} r)$, where $\epsilon, r, q \in \mathbb{Q}_{[0,1]}$, $r > \epsilon$, $q < \delta_{f,i}(\epsilon)$,

$mp$: $\dfrac{\varphi \quad \psi \mathbin{\dot{-}} \varphi}{\psi}$,

$gen$: $\dfrac{\varphi}{\sup_x \varphi}$.

## 2.3 Reducing CL to RPL∀

We can consider **CL** as a theory in **RPL∀**. Every connective in **CL** can be defined in **RPL∀**, $\varphi \mathbin{\dot{-}} \psi := \psi \to \varphi$, $\sup_x \varphi := \forall x \varphi$, $\bar{q} := \bar{q}$. Axioms and rules are identical except for $S4$ and $UL$. We assume that $S4$ is not part of **LL∀** or **RPL∀**, unless explicitly stated otherwise. We can put all $UL$ axioms in our theory $T$, instead of adding it as an axiom for **RPL∀**. Note that $UL$ axioms are actually expressed in the language of **RPL∀**, and also form a decidable subset of its sentences. Given a theory $T$ in **CL** we build $T^*$ by adding all $UL$ axioms to $T$, which is an effective process. Any model of $T^*$ in **RPL∀** will be a model of $T$ in **CL** since in any model of this theory, satisfaction of $UL$ forces every function and relation to be uniformly continuous with respect to their given modulus of continuity. Therefore a problem in **CL** is reduced to a problem in **RPL∀**.

## 2.4 Examples of Continuous Structures and Theories

We give a few metric structures and continuous theories as examples from [BBHU]:

1. A complete, bounded metric space $(M, d)$ with no additional structure.

2. A structure $M$ in the usual sense from first-order logic. $d$ is the discrete metric on the underlying set, i.e. $d(a,b) = \begin{cases} 0 & a=b \\ 1 & a \neq b \end{cases}$ and a relation is considered as a predicate taking values (truth-values) in the set $\{0,1\}$. So, **CL** is, in a sense, a generalization of first-order logic.

3. Let $C$ be the set of continuous functions over $[0,1]$ and $\int$ as a unary predicate, defined as $f \mapsto \int_0^1 f(x)dx$. Then $(C, \int, ||.||_\infty)$ is a metric structure.

4. $(C, <\cdot,\cdot>, ||.||_\infty)$, where $C$ is the same as above, and $<\cdot,\cdot>$ is the $L_2$ inner product over $C$, i.e. $<f,g> := \int_0^1 f(x)g(x)dx$.

5. An infinite dimensional Hilbert space over $\mathbb{R}$ can be considered as a metric structure, according to [BBHU, §15], if only the range of its metric is assumed to be bounded. This obstacle can be overcome easily by introducing many-sorted models. We will not study many-sorted structures here, although generalization of our results is straightforward. Let $H$ be an infinite dimensional Hilbert space over $\mathbb{R}$, i.e. $H$ is a complete normed space with an inner product $<\cdot,\cdot>$ which induces the norm. $H$ can be considered as a metric structure $(B_n : n \geq 1, I_{m,n} : m < n, 0, \lambda : \lambda \in \mathbb{R}, +, -, <\cdot,\cdot>, ||.||)$ where $B_n := \{x \in H : ||x|| \leq n^2\}$, $I_{m,n} : B_m \hookrightarrow B_n$ is the inclusion map for $m < n$. It is proved in [BBHU, §15] that $H$ has quantifier elimination property and for any $H_1$ and $H_2$, $Th(H_1) = Th(H_2)$, therefore $H$ is semantically linear complete.

6. The unit ball $B$ of a Banach space $X$ over $\mathbb{R}$ or $\mathbb{C}$: as functions we may take the maps $f_{\alpha\beta}$, defined by $f_{\alpha\beta}(x,y) = \alpha x + \beta y$, for each pair of scalars satisfying $|\alpha| + |\beta| \leq 1$; the norm may be included as a predicate, and we may include the additive identity $0_X$ as a distinguished element. See [BU, 4.4, 4.5].

7. Probability spaces can be represented as metric structures using their measure algebra [BBHU, §16]. Assume $(X, B, \mu)$ is a probability space, i.e. $X$ is an arbitrary set, $B$ is a $\sigma$-algebra of subsets of $X$, and $\mu$ is a $\sigma$-additive measure in $B$ such that $\mu(X) = 1$. $a \in B$ is an *atom* iff there is no $b \in B$ such that $b \subseteq a$ and $0 < \mu(b) < \mu(a)$. $a \in B$ is atomless iff no subset of it is an atom. $(X, B, \mu)$ is atomless iff $X$ is atomless. For any atomless $a \in B$ and $r \in [0,1]$ there exists $b \in B$ such that $b \subseteq a$ and $\mu(b) = r \cdot \mu(a)$. $a, b \in B$ represent the same *event* iff $\mu(a\Delta b) = 0$ ($\Delta$ being symmetric difference, i.e. $a\Delta b := (a \cap b^c) \cup (a^c \cap b)$), in which case we say $a \sim_\mu b$. $\sim_\mu$ is an equivalence relation. Denote equivalence class of $a$ by $[a]_\mu$, and $\hat{B}$ as set of equivalence classes, which is called *measure algebra* of $(X, B, \mu)$. Since $\sim_\mu$ is congruence with respect to complement, union, intersection, and $\mu$, these operations are well-defined on $\hat{B}$. We also have that $\hat{B}$ is a $\sigma$-algebra and $\mu$ is a strictly positive countably additive measure over it. The metric structure $M = (\hat{B}, 0, 1, \cdot^c, \cup, \cap, \mu, d)$ is called a probability structure, where $d([a]_\mu, [b]_\mu) := \mu(a\Delta b)$, $0 := [\emptyset]_\mu$, and $1 := [X]_\mu$. Modulus of uniform continuity of $\cdot^c$, $\cup$, and $\cap$ are respectively $\epsilon \mapsto \epsilon$, $\epsilon \mapsto \epsilon/2$, and $\epsilon \mapsto \epsilon/2$.

It turns out that theory of a probability structure is axiomatizable:

(a) Boolean Algebra Axioms: axioms are closure of equations, thus expressible in **CL**

(b) Measure Axioms:

$\mu(0) = 0$: $\mu(0)$,
$\mu(1) = 1$: $\neg\mu(1)$,
$\sup_{x,y}(\mu(x \cap y) \dotminus \mu(x)) = 0$: $\forall x, y \, (\mu(x) \to \mu(x \cap y))$,
$\sup_{x,y}(\mu(x) \dotminus \mu(x \cup y)) = 0$: $\forall x, y \, (\mu(x \cup y) \to \mu(x))$,
$\sup_{x,y} |(\mu(x) \dotminus \mu(x \cap y)) - (\mu(x \cup y) \dotminus \mu(y))| = 0$: $\forall x, y \, ((\mu(x \cap y) \to \mu(x)) \leftrightarrow (\mu(y) \to \mu(x \cup y)))$,

(c) Connection between $\mu$ and $d$: $\sup_{x,y} |d(x, y) - \mu(x\Delta y)| = 0$,
i.e. $\forall x, y \, (x \approx y \leftrightarrow \mu(x\Delta y))$,

(d) Atomlessness: $\sup_x \inf_y |\mu(x \cap y) - \mu(x \cap y^c)| = 0$,
i.e. $\forall x \, \exists y \, (\mu(x \cap y) \leftrightarrow \mu(x \cap y^c))$.

The last three axioms of Measure Axioms express $\mu(x) + \mu(y) = \mu(x \cup y) + \mu(x \cap y)$. The theory $PR_0$ satisfies (a), (b), and (c), and $PR$ satisfies (a), (b), (c), and (d). For every metric structure $M$ in this language we have that $M \vDash PR_0$ iff $M$ is the probability structure of some probability space $(X, B, \mu)$. $PR$ has quantifier elimination, and is separably categorical, i.e. every two separable metric models of it are isomorphic, and by using downward Löwenheim-Skolem theorem for **CL** one can obtain a separable elementary substructure of any given structure, therefore, $PR$ is semantically linear. We will show that every axiomatizable semantically linear theory has a computable model, thus $PR$ has a unique computable model.

8. The $L^p$ Banach lattice is another example of [BBHU, §17], [BBH]. Similarly to PR, it is axiomatizable, has quantifier elimination, and is separably categorical, therefore is semantically linear complete, and has a unique computable model.

## 2.5 Computable Analysis and Continuous Domain Theory

We will closely follow [EH98, ES99b][7]. We do not use any special property of effective continuous domain approach to computable analysis in our work, and therefore transferring our results to other approaches, such as TTE [Wei00] or Scott-Ershov $\omega$-algebraic domains [SHT07], is straightforward. We prefer continuous domains since they are simpler and more intuitive.

---

[7]See also [AJ94]

We will use fairly standard notions from computability theory. $\mathbb{N}$ denotes nonnegative integers. $\varphi_n$ is the $n$th partial computable function, and $W_n := dom(\varphi_n)$, the $n$th computably enumerable (c.e.) subset of $\mathbb{N}$. $<\cdot,\cdot>$, $\pi_0$, and $\pi_1$ are the usual Cantor's pairing function and its inverses, satisfying $<\pi_0(n), \pi_1(n)> = n$, $\pi_i(<n_0, n_1>) = n_i$ for $i = 0, 1$. Any c.e. $A$ can be expressed as a computable union of decidable sets $\{A_n\}_{n \in \omega}$. $K$ is the associated set with halting problem, i.e., $K = \{n : \varphi_n(n) \text{ is defined}\}$.

Let $(P, \sqsubseteq)$ be a partially ordered set (poset), i.e. $\sqsubseteq$ is reflexive, antisymmetric, and transitive. We think of $P$ as partial information and $\sqsubseteq$ as refinement (having more information). Upper and lower bounds, maximal, minimal, top (maximum, $\top$), bottom (minimum, $\bot$), supremum ($\sqcup$) and infimum ($\sqcap$) with respect to $\sqsubseteq$ are defined as usual. The set of maximal elements of $P$ is denoted by $\max(P)$. Assume that $A$ is a subset of $P$. We say that $A$ is a *chain*, if every two elements of $A$ are comparable; and $A$ is *directed*, if every two elements have an upper bound in $A$. *Upper set* of $A$ ($\uparrow A$) is $\{x \in P : \exists y \in A \; y \sqsubseteq x\}$. *Lower set* of $A$ ($\downarrow A$) is defined similarly.

A partial order is called *dcpo* (directed complete partial order), or *bcpo* (bounded complete partial order), if it contains sups of all directed sets, or if it contains all bounded sets, respectively. We call a dcpo *pointed* if it contains a least element $\bot$. A *cpo* (complete partial order) is a pointed dcpo.

If $S$ is a set, the *flat domain of* $S$, $S_\bot$, is the partial order over $S \cup \{\bot\}$, where $\bot$ is minimum and no two elements of $S$ are comparable. For example, if $B$ is the bi-valued Boolean set $\{0, 1\}$, $B_\bot$ is flat domain of $B$.

We say that $x$ is *way-below* $y$, or $x$ *approximates* $y$, and write $x \ll y$, iff for all directed $A$, $y \sqsubseteq \sqcup A$ implies $x \sqsubseteq a$ for some $a \in A$. An element which approximates itself is called *compact* or *finite*. We denote the set of all elements way-below $x$ (way-above $x$) by $\Downarrow x$ ($\Uparrow x$). A *base* for $D$ is a subset $B$ of $D$ such that for all $x \in D$ we can approximate $x$ by elements of $D$, i.e. $\Downarrow x \cap B$ is directed and $x = \sqcup(\Downarrow x \cap B)$. A dcpo is called ($\omega$-)*continuous* iff it has a (countable) base. A continuous dcpo is called a *continuous domain* or simply a domain.

Real interval domain, $I\mathbb{R}$, is the set of all closed intervals in $\mathbb{R}$ including $\mathbb{R}$ itself, partially ordered by reverse inclusion relation. Unit interval domain, $I[0, 1]$, is defined similarly. It is easy to see that both $I\mathbb{R}$ and $I[0, 1]$ are $\omega$-continuous domains, and $x \ll y$ iff $x^o \supseteq y$, where $x^o$ is the interior of $x$. The set of maximal elements are homeomorphic to the original space by $x \mapsto \{x\}$.

Assume that $(D, \sqsubseteq)$ and $(E, \sqsubseteq)$ are two dcpos. We say that a function $f : D \to E$ is *Scott-continuous* (or simply continuous), iff:

- monotonicity: for all $x, y \in D$, $x \sqsubseteq y \to f(x) \sqsubseteq f(y)$

- continuity: for all directed subsets $A$ of $D$, $f(\sqcup A) = \sqcup f(A)$

We denote the set of all continuous functions from $D$ to $E$ by $[D \to E]$, which is itself a dcpo by pointwise ordering, i.e. $f \sqsubseteq g$ iff $\forall x \in D \; f(x) \sqsubseteq g(x)$. $f \in [D \to E]$ is *strict* iff $f(\bot_D) = \bot_E$.

The product of two domains is also a domain, where $(x_0, y_0) \sqsubseteq (x_1, y_1)$ iff $x_0 \sqsubseteq x_1 \wedge y_0 \sqsubseteq y_1$. It is also easy to check that $(x_0, y_0) \ll (x_1, y_1)$ iff $x_0 \ll x_1 \wedge y_0 \ll y_1$; a function from product is continuous if it is continuous in each variable separately.

Let $D$ be a $\omega$-continuous dcpo with a given countable base $B = \{b_n : n \in \omega\}$. Without loss of generality, assume that $b_0 = \bot$. We say that $D$ is *computably given with respect to* $b$ iff $\{(m, n) : b_m \ll b_n\}$ is a c.e. set. An element $x \in D$ is *computable* iff $\{n : b_n \ll x\}$ is computably enumerable, i.e. $\{n : b_n \ll x\}$ is the range of a computable function; thus $x$ will be supremum of this c.e. set. A sequence of elements of $D$, $\{x_n : n \in \omega\}$ is *computable* iff there exists a computable function $\varphi : \omega \times \omega \to \omega$ such that $\{m : b_m \ll x_n\} = \{\varphi(n, k) : k \in \omega\}$, i.e. approximations of the elements of the sequence are uniformly computable. If $D$ and $E$ are computable domains, then their product is a computable domain too.

Let $f \in [D \to E]$, where $D$ and $E$ are computably given by bases $A = \{a_n : n \in \omega\}$ and $B = \{b_n : n \in \omega\}$. We say that $f$ is *computable* iff $f$ takes computable sequences in $D$ to computable sequences in $E$, i.e. $\{(m, n) : b_m \ll f(a_n)\}$ is computably enumerable.

Let $I[0, 1] = \{[p, q] : p, q \in [0, 1]\}$ where $[p, q] \sqsubseteq [r, s]$ iff $[p, q] \supseteq [r, s]$. Elements of $I[0, 1]$ are approximations to real numbers in the unit interval. The order of approximation is given by $[p, q] \ll [r, s]$ iff $[p, q] \supseteq (r, s)$. As an effective base, we use $B = \{[p, q] : p, q \in \mathbb{Q}_{[0,1]}\}$. The bottom element of this ordering is $\bot = [0, 1]$. For $r \in \mathbb{Q}_{[0,1]}$, $\{r\}$ is a maximal element. All maximal elements are of this form. The unit interval is embedded in $I[0, 1]$ by $\{\cdot\} : [0, 1] \to I[0, 1]$, taking $x \mapsto \{x\}$.

Let $(X, d)$ be a metric space. *Formal ball* construction on $X$ is defined as $BX := \{(x, r) : x \in X, r \in \mathbb{R}^*\}$, where $\mathbb{R}^*$ is the set of nonnegative real numbers. We define $(x, r) \sqsubseteq (y, s)$ iff $d(x, y) \leq r - s$. $(BX, \sqsubseteq)$ is a poset. Many notions for dcpo can be generalized for posets, e.g. continuity, by replacing quantifications over "directed" with "directed with an upper bound." We now state some important properties of Formal Balls: $BX$ is always a continuous dcpo. $(X, d)$ is a complete metric space iff $BX$ is a dcpo. $(x, r) \ll (y, s)$ iff $d(x, y) < r - s$. If $A$ is a dense subset of $X$ and $Q$ is a dense subset of $\mathbb{R}^*$, then $A \times Q$ is a base for $BX$. $X$ is separable iff $BX$ is $\omega$-continuous. $X$ is embedded in $BX$ by $x \mapsto (x, 0)$, where each $(x, 0)$ is a maximal element of $BX$. This embedding is continuous with respect to Scott and Lawson topologies.

Therefore, if $X$ is a separable complete metric space, then $BX$ is $\omega$-continuous dcpo. $X$ is a computable metric space if its domain is computable. Maximal elements will be homeomorphic to $X$. We denote this homeomorphism by $\iota_X : Max(D) \to X$. We adopt some concepts from [SHT07] to apply to continuous domains. Let $X$ and $Y$ be metric spaces represented by $D$ and $E$. A function $f : X \to Y$ is represented by $\hat{f} : D \to E$ iff $\hat{f}[Max(D)] \subseteq Max(E)$ and for any $x \in Max(D)$, $\iota_Y(\hat{f}(x)) = f(\iota_X(x))$. $f$ is a computable function iff $\hat{f}$ is a computable function between corresponding domains.

$(X_0, X_1, ..., f_0, ...)$ is called a many-sorted *topological structure* when $X_i$s are topological spaces and $f_j$s are continuous functions over them. A (pre-)metric structure is a special case, where $\mathbb{R}$ is one of the spaces, and for every space $X$ in the structure we have a continuous function $d_X$ that is a (pre-)metric on $X$. Definition of a domain presentation is straightforward. A structure is computably presentable iff there is a computable domain presentation for it.

## 2.6 Effective Model Theory

A countable, first order language $L$ is said to be *computably given* iff the sets of variables, functions, relations are decidable and arity of functions and relations are computable. Let $T$ be a theory in $L$. $T$ is *computably axiomatizable* or simply *axiomatizable* iff the set of axioms are computably enumerable. $T$ is *decidable* iff the set of theorems of $T$ is decidable.

A countable structure is called *computably representable* iff its relations are decidable and its functions are computable, i.e. its *atomic diagram*, $Diam_{at}(M) := \{\varphi \in L_M : M \vDash \varphi \land \varphi \text{ is atomic}\}$, is decidable. A countable structure is called decidable iff it is computably presented and satisfaction is decidable. Equivalently, its *Elementary Diagram*, $Diag_{el}(M) := \{\varphi \in L_M : M \vDash \varphi\}$, is decidable. When adopting this notion to logics with an uncountable set of truth values, we need to clarify what we mean by decidable and computable. We use effective domain theory. Decidability of relations in classical models is equivalent to computability of their characteristic functions. Extending this to continuous case, where the range of characteristic functions is $[0,1]$, in a computable continuous structure, we require characteristic functions of relations to be computable. Therefore, a computable continuous model of a theory is a computable continuous structure that is also a model of the theory.

For countable models, accordingly, equivalent formulations by $Diag_{at}(M)$ and $Diag_{el}(M)$ can be adopted by replacing classical decidability of truth for formulas with computability of their truth degree. However, we want to allow a complete metric space as the underlying set of models of continuous logic, and since most of these structures are uncountable, we cannot use this formulations directly. For complete metric spaces with continuous operations and relations, we can use a countable dense set for approximating elements, and consider only formulas that use parameters only from this dense subset. Since operations and relations are continuous, every formula in the full diagrams can be approximated with arbitrary precision. Assume that $M$ is a complete metric space. We say that $M$ is computably representable iff there is a countable dense subset $X$ of $M$ so that $Diag_{at}(M, X) := \{(\varphi, x) \in L_X \times B : x \ll ||\varphi||_M \land \varphi \text{ is atomic}\}$ is computably enumerable, where $B$ is the effective base of the unit interval discussed in previous section. $M$ is decidable iff $Diag_{el}(M, X) := \{(\varphi, x) \in L_X \times B : x \ll ||\varphi||_M\}$ is computably enumerable. (The construction of a computably presented structure from a computable enumeration of elementary diagram is straightforward.)

Note that every linear complete theory that has a decidable model is weakly axiomatizable.

# 3  Constructive Henkin Construction

In this section, we will construct an effective model for an effectively given theory.

**Theorem 3.1.** *Every consistent, complete, computably axiomatizable theory in $\mathbf{RPL}\forall$ is decidable.*

*Proof.* Assume $T$ is a theory as stated. Because $\mathbf{RPL}\forall$ is an effective deduction system with computably enumerable axioms and two deduction rules, we can effectively list all valid proofs with assumptions in $T$. Since $T$ is complete and consistent, we have $T \nvdash \varphi \Leftrightarrow T \vdash \neg\varphi$. For any given sentence such as $\varphi$, either $T \vdash \varphi$ or $T \vdash \neg\varphi$ has a proof, and, at some point, that proof will appear in the list. According to whether it is a proof of $\varphi$ or $\neg\varphi$, we answer positively or negatively. $\square$

Completeness is too strong to be satisfied by a theory. The Lindenbaum Algebra of a complete theory will have just two classes, in other words, the logic will be the classical bi-valued logic. We want to allow truth degrees in $[0,1]$, not just in $\{0,1\}$. Linear-completeness is an acceptable condition, which forces the classes of the Lindenbaum Algebra of a theory to be linearly ordered. The concepts of *decidable model* and *decidable theory* are not as useful for many-valued logics as for the classical bi-valued one. Decidability is the same as saying the *characteristic functions* of truth and validity are computable. The natural extensions of these concepts to many-valued cases are the *computability of the degree* of the truth and provability of formulas. Crisp truth and validity are not computable functions because of discontinuity ([Wei00, §2.4]). It is proved in [Haj98, Theorem 6.3.15] that crisp truth in $\mathbf{LL}\forall$ is $\Pi_2-complete$.

**Theorem 3.2.** *For every consistent, linear-complete, computably axiomatizable theory in $\mathbf{RPL}\forall$, the provability degrees of sentences is computably comparable.*

*Proof.* The proof is similar to previous one. Since $T$ is linear-complete, for any given two sentences such as $\varphi$ and $\psi$, at least one of $T \vdash \varphi \to \psi$ and $T \vdash \psi \to \varphi$ has a proof, and, at some point, that proof will appear in the list. According to whether it is a proof of $T \vdash \varphi \to \psi$ or $T \vdash \psi \to \varphi$, we answer $|\psi|_T \leq |\varphi|_T$, or $|\varphi|_T \leq |\psi|_T$. $\square$

We use Henkin's technique for model construction. To construct a decidable model, we need to deal effectively with real numbers. For this purpose, we need a model for computable analysis.

The provability degree of a formula with respect to a theory is a real number, and since comparison of real numbers is not decidable in finite time, we have to accept existence of distinct elements in our model which are equal with provability degree 0. The underlying set for a computably representable (which is the weakest computability notion) model in classical computable model theory is natural numbers (although often not explicitly stated) and, therefore, has a decidable equality.

**Theorem 3.3.** *There is a computably axiomatizable, consistent theory that does not have a computable model (with a decidable equality on its domain), where equality is interpreted as identity.*

*Proof.* Assume that $a$ is a computable real number in $[0,1]$ such that '$a \stackrel{?}{=} 0$' is not decidable.[8] Suppose $\{s_i\}_{i\in\omega}$ ($\{r_i\}_{i\in\omega}$) is a decreasing (increasing) computable sequence of rationals converging to $a$. Let $c$ and $d$ be two distinct constants of the language, and $T = \{\overline{s_i} \to c \approx d : i \in \omega \land s_i \in \mathbb{Q}_{[0,1]}\} \cup \{c \approx d \to \overline{r_i} : i \in \omega \land r_i \in \mathbb{Q}_{[0,1]}\}$. $T$ is computably axiomatized. If there existed a computable model $M \vDash T$, then by checking whether the interpretation of $c$ and $d$ are identical, we could decide '$a \stackrel{?}{=} 0$;' a contradiction. □

This can be avoided at the expense of allowing computable models over domains that lack a decidable equality. For countable domains, this is very artificial. But computable uncountable domains essentially lack a decidable equality, since they are represented by a countable sequence of approximating elements from a countable domain.

Before proceeding, we prove a lemma.

**Lemma 3.4.** *Provability degree of a formula with respect to a consistent, linear-complete, computably axiomatizable (Henkin) theory is computable.*

*Proof.* Since $T$ is linear-complete, $|\varphi|_T = \sqcap\{r \in \mathbb{Q}_{[0,1]} : T \vdash \overline{r} \to \varphi\} = \sqcup\{r \in \mathbb{Q}_{[0,1]} : T \vdash \varphi \to \overline{r}\}$, which is a maximal element of $I[0,1]$, using the embedding $x \mapsto \{x\}$.

Let $B = \{[r,s] \in I[0,1] : r < s \in \mathbb{Q}_{[0,1]}\}$. $B$ is an effective base for the unit interval domain $I[0,1]$. We have $|\varphi|_T = \bigsqcup_{I[0,1]}\{[r,s] \in B : (T \vdash \overline{s} \to \varphi) \land (T \vdash \varphi \to \overline{r})\}$. It is easy to check that this set is directed, and effectively listable. Let $\{q_i\}_{i\in\omega}$ be an effective list of rationals in $[0,1]$. The following algorithm gives an enumeration of the set:

1. list $[0,1]$;

2. $r \leftarrow 0; s \leftarrow 1$;

3. for $i \in \omega$,
   if $r \leq q_i \leq s$, depending on whether $T \vdash \varphi \to \overline{q_i}$ or $T \vdash \overline{q_i} \to \varphi$, go to (a) or (b)

   (a) $s \leftarrow q_i$; for $j < i$ and $q_j < s$, list $[q_j, s]$;
   (b) $r \leftarrow q_i$; for $j < i$ and $r < q_j$, list $[r, q_j]$;

For any $r \in \mathbb{Q}_{[0,1]}$, either $T \vdash \varphi \to \overline{r}$ or $T \vdash \overline{r} \to \varphi$. Listing the proofs, we will find, at some point, which is the case. Therefore, $|\varphi|_T = \bigsqcup_{I[0,1]}\{[r,s] \in B : (T \vdash \overline{s} \to \varphi) \land (T \vdash \varphi \to \overline{r})\}$ is computable, using the effective base $B$ of $I[0,1]$ as defined above. □

**Theorem 3.5.** *Every consistent, linear-complete, computably axiomatizable Henkin theory in* **RPL∀** *has a decidable model.*

---

[8]See [TvD87], [Bee80], and [Wei00].

*Proof.* Assume that $T$ is as stated, and $C = \{c_i\}_{i \in \omega}$ is the set of constants of the language of $T$. First, note that all of the logical connectives, specially $\rightarrow$, commute with provability degree, and for a Henkin theory, we have $|\forall x\, \varphi|_T = \sqcup\{|\varphi[c/x]|_T : c \in C\}$.

Let the underlying set of our model, $M$, be $C$. Assuming the language of $T$ is decidable, $M$ is a decidable set. As we stated before, we cannot expect our model to identify elements that are equal with provability degree 0. The notion of *similarity*,[9] i.e. fuzzy equality, will be interpreted by a pseudo-metric $\rho$ on $M$. We define $\rho(x, y) = |x \approx y|_T$. Similarity is computable since provability degree with respect to $T$ is computable. It is easy to check that $\rho$ is a pseudo-metric, because of similarity axioms (2.1).

$$\frac{\overline{\phantom{xx}}\ S1}{\frac{x \approx x}{\overline{0} \rightarrow x \approx x}} \qquad \frac{\overline{\phantom{xx}}\ S2}{\frac{\overline{r} \rightarrow x \approx y,\ x \approx y \rightarrow y \approx x}{\overline{r} \rightarrow y \approx x}}$$

$$\frac{\frac{\overline{r_1} \rightarrow x \approx y, \overline{r_2} \rightarrow y \approx z}{\overline{r_1} \odot \overline{r_2} \rightarrow x \approx y \odot y \approx z}\ ,\ \overline{\phantom{xxxxxxxxx}}\ S3}{\overline{r_1 \dot{+} r_2} \rightarrow x \approx z}$$

Hence, $\rho(a, a) = 0$, $\rho(a, b) \leq \rho(b, a)$, $\rho(a, c) \leq \rho(a, b) \dot{+} \rho(b, c) \leq \rho(a, b) + \rho(b, c)$. For each function symbol $f$ of $T$, we have $T \vdash \forall \overrightarrow{x}\, \exists y\, f(\overrightarrow{x}) \approx y$, and therefore for any $\overrightarrow{a} \in M$, $T \vdash \exists y\, f(\overrightarrow{a}) \approx y$. Since $T$ is Henkin, there is a constant $c$ in $T$ so that $T \vdash f(\overrightarrow{a}) \approx c$. We define the interpretation of $f$ over $\overrightarrow{a}$ in $M$, $f^M(\overrightarrow{a})$, as the first constant $c$ where we find a proof of $T \vdash f(\overrightarrow{a}) = c$. This gives a computable interpretation of $f$ in $M$. This procedure can be used repeatedly for computing complex terms.

Interpret $\overline{r}^M$ by $r$, for $r \in \mathbb{Q}_{[0,1]}$. For any relation symbol $R$ of $T$, define $R^M(a)$ as $|R(a)|_T$, which is again computable. Actually, we can compute the value of any sentence of $M$ directly using provability degree.

So far, we have created a decidable structure $M$. Using commutativity of provability degree with logical connectives, we now show that $M$ is indeed a model of $T$, i.e. $M \vDash T$. We prove the stronger $|\varphi|_T = ||\varphi||_M$ by induction on formulas.

For any $r, s \in \mathbb{Q}_{[0,1]}$, by axiom $R$ of 2.1, $T \vdash \overline{r} \rightarrow \overline{s}$ iff $T \vdash \overline{s \dot{-} r}$; by strong consistency (which is equivalent to consistency), that is equivalent to $s \dot{-} r = 0$, i.e. $s \leq r$.

Let $q \in \mathbb{Q}_{[0,1]}$, then

$$|\overline{q}|_T = \bigsqcup_{I[0,1]} \{[r, s] \in B : T \vdash \overline{q} \rightarrow \overline{r} \wedge T \vdash \overline{s} \rightarrow \overline{q}\} = \bigsqcup_{I[0,1]} \{[r, s] \in B : r \leq q \leq s\} = \{q\}$$

---

[9] Actually we use a pseudo-metric, because of the isomorphism $x \mapsto 1 - x$. This transformation maps an interpretation of similarity to a pseudo-metric and vice versa; see [Haj98, §5.6].

On the other hand, $||\bar{q}||_M$ is defined as $q$.

For relational symbol $R$, we have $||R(\vec{a})||_M = |R(\vec{a})|_T$, by definition. Similarly, for $\approx$.

For $\varphi = \varphi_1 \to \varphi_2$,

$$\begin{aligned}
||\varphi||_M &= \\
||\varphi_1 \to \varphi_2||_M &= & \text{(by definition)} \\
||\varphi_2||_M \mathbin{\dot{-}} ||\varphi_1||_M &= & \text{(by induction hypothesis)} \\
|\varphi_2|_T \mathbin{\dot{-}} |\varphi_1|_T &= & \text{(by lemma 2.6)} \\
|\varphi_1 \to \varphi_2|_T &= \\
|\varphi|_T &
\end{aligned}$$

For $\varphi = \forall x\, \psi$,

$$\begin{aligned}
||\varphi||_M &= \\
||\forall x\, \psi||_M &= & \text{(by definition)} \\
\sqcup_{a \in M} ||\psi[x/a]||_M &= & \text{(by induction hypothesis)} \\
\sqcup_{a \in M} |\psi[x/a]|_T &= & (M = C) \\
\sqcup_{a \in C} |\psi[x/a]|_T &= & \text{(by lemma 2.6)} \\
|\forall x\, \psi|_T &= \\
|\varphi|_T &
\end{aligned}$$

Therefore the above interpretations over $M$ give a decidable model of $T$. □

**Remark 3.6.** *If the theory contains S4, the resulting model will be extensional.*

We state two results about completion process of a theory.

**Theorem 3.7.** *Every computably axiomatizable, consistent theory in* **RPL∀** *has a K-computably axiomatizable, consistent, linear-complete Henkin extension.*

*Proof.* Assume that $T$ is as stated. $L^*$ is created by adding $\{c_n\}_{n \in \omega}$ to $L$. Let $\{(\varphi_n, \psi_n)\}_{n \in \omega}$ be an enumeration of pairs of $L^*$-sentences, so that any sentence $\varphi$ appears in $\{\varphi_n\}_{n \in \omega}$ infinitely many times. We build $T^*$ in $\omega$ steps:

1. Let $T_0 \leftarrow T$;

2. For $n \in \omega$,

(a) If $T_n \nvdash \varphi_n \to \psi_n$, then let $T_{n+1} \leftarrow T_n \cup \{\psi_n \to \varphi_n\}$;

(b) If $\varphi_n$ is of the form $\exists x\, \psi$, and $T_n \vdash \exists x\, \psi$,
take an unused constant $c$, and put $T_{n+1} \leftarrow T_n \cup \{\psi[x/c]\}$;

3. Let $T^* = \cup_{n \in \omega} T_n$.

$T^*$ is linear complete, since for any pair of sentences such as $(\varphi_n, \psi_n)$, either $T_n \vdash \varphi_n \to \psi_n$ and thus $T^* \vdash \varphi_n \to \psi_n$, or $T_n \nvdash \varphi_n \to \psi_n$, and thus $\psi_n \to \varphi_n \in T_{n+1} \subseteq T^*$ and $T^* \vdash \psi_n \to \varphi_n$.

To see $T^*$ is Henkin, assume that $T^* \vdash \exists x\, \psi$. There is an $m \in \omega$ such that $T_m \vdash \exists x\, \psi$. $\varphi_n = \exists x\, \psi$ for some $m \leq n$. $T_n \vdash \exists x\, \psi$, therefore $\psi[x/c] \in T_{n+1} \subseteq T^*$ for some $c$, and $T^* \vdash \psi[x/c]$.

To see $T^*$ is consistent, it is sufficient to prove that if $T_n$ is consistent, then $T_{n+1}$ is consistent. $T_0$ is consistent. Assume that $S \cup \{\varphi \to \psi\} \vdash \overline{1}$, i.e. is inconsistent. By 2.1, there is an $n \in \omega$ such that $S \vdash (\varphi \to \psi)^n \to \overline{1}$. By 2.7, $S \vdash (\varphi \to \psi)^n \vee (\psi \to \varphi)^n$, therefore $S \vdash \overline{1} \vee (\psi \to \varphi)^n$, $S \vdash (\psi \to \varphi)^n$, $S \vdash (\psi \to \varphi)$.

On the other hand, if $S \vdash \exists x\, \psi$, $S \nvdash \overline{1}$, and $c$ does not occur in $S$, then $S \cup \{\psi[x/c]\}$ must be consistent.

The above algorithm is computable in $K$, thus it is in class $\Delta_2^0$, i.e. it is *limit computable*. □

**Remark 3.8.** *It is easy to check that if the theory were to be computably enumerable in $A$, the procedure above would yield a theory computably enumerable in $A'$, and, therefore, would have a $A'$-decidable model.*

We now show that this is the best result one can achieve. We actually prove a much stronger result by showing that for theories in classical logic, i.e. those containing $\{\varphi \vee \neg\varphi : \varphi \in L\}$, a classical decidable model can effectively be generated from our construction. This shows that computable continuous model theory is a generalization of classical computable model theory, as continuous logic is a generalization of classical logic.

**Theorem 3.9.** *There is an algorithm that given a classical theory $T$ in $\mathbf{RPL}\forall$, and a countable decidable model $M \vDash T$, produces a classical computable model.*

*Proof.* Take $M$ as our underlying set of the model. Since the theory is classical, $M$ is a classical model (just forget that the range is $[0,1]$, since all values fall in $\{0,1\}$). It is sufficient to show that its elementary diagram is decidable. Let $\varphi \in L_M$ be a sentence. Since $T$ is classical, it proves PEM for all instances of $\varphi$ and we have $T \vdash \varphi \vee \neg\varphi$, therefore $M \vDash \varphi \vee \neg\varphi$, i.e. $\min\{||\varphi||_M, 1 - ||\varphi||_M\} = 0$, or $||\varphi||_M$ is either 0 or 1. Since $||\varphi||_M$ is a maximal element, at some point of $||\varphi||_M$'s computation, we will see that the approximation does not contain either 0 or 1. Since $||\varphi||_M$ is either 0 or 1, we can say at this point which of them it is, and can decide whether it is true or false in our model. □

It is now easy to get the following corollary:

**Corollary 3.10.** *There exists a computably axiomatizable, consistent theory in $\mathbf{RPL}\forall$ (and thus in $\mathbf{LL}\forall$) that has no decidable model.*

*Proof.* There is a consistent theory in classical logic which is computably axiomatizable but does not have a decidable classical model by [Mil99]. Since $\mathbf{RPL}\forall$ is a conservative extension of classical logic, this theory is consistent in $\mathbf{RPL}\forall$. By the previous theorem, it can't have a decidable continuous model. [10] □

**Corollary 3.11.** *The theory discussed above has no computably axiomatizable, consistent, linear-complete extension.*

*Proof.* Take the theory in the previous theorem. If it had a computably axiomatizable consistent linear-complete extension, then by 3.5, it would have a decidable model, contradicting 3.10. □

This technique can be used easily to carry over negative results from computable classical model theory to computable continuous model theory.

Our last theorem is a construction of a decidable complete metric model from the countable metric model produced by a Henkin construction.

**Theorem 3.12.** *Given a countable, decidable model of $T$ as in the 3.5, there is an elementarily equivalent, decidable, complete, metric model of $T$.*[11]

*Proof.* Take $M$ as the metric completion of the underlying countable pseudo-metric space $B$ of our Henkin model after taking the quotient with respect to its pseudo-metric. Functions and relations of $B$ can be extended to $M$ by taking limits of their values over elements of $B$ (a countable dense subset of $M$) since they are continuous. An easy induction shows that $B$ is an elementary substructure of $M$. Therefore, for every $\varphi \in L_B$, $B \vDash \varphi$ iff $M \vDash \varphi$. Thus $Diag_{el}(M,B) = Diag_{el}(B)$, and is computably enumerable. □

**Remark 3.13.** We have shown in Theorem 3.5 that any computable linear complete theory has a decidable model. It is easy to see the essential property of a linear theory we use is the computability of the truth degree of the formulas. Therefore with minor adjustments, the proofs will hold also for semantically linear complete theories in continuous logic. Hence the following corollary follows.

**Corollary 3.14.** *The unique separable models of probability spaces and $L^p$ Banach lattices are decidable.*

---

[10] The idea of proof is the well-known fact in constructive mathematics that $(\forall x < y \in \mathbb{R})(\forall z \in \mathbb{R}), (x < z \lor z < y)$. [TvD87].

[11] This is computable version of one of fundamental results of continuous logic, e.g., see [BU, 2.10].

# 4 Conclusions and further works

In this paper we defined basic notions of computability for continuous model theory. We proved a constructive version of Henkin's construction. A paper in preparation will discuss constructive versions of *omitting types theorem, elementary chain construction* , and *existence of prime and saturated models.* The work on other theorems and concepts such as *computably categorical* theories is still in progress.

Extending definitions and theorems to multi-sorted case is straightforward. Another interesting topic is study of partial models. Partial models for **CL** and **LL**$\forall$, are what Boolean valued models are for classical first order logic. The capturing of a satisfactory definition for a partial model is in a sense more complicated than that of the definition of Boolean valued models. The basic trick is to require $||\varphi||_M \dotdiv ||\psi||_M \sqsubseteq ||\varphi \dotdiv \psi||_M$ in the domain theoretic sense. For example, $[0,1] \dotdiv [0,1] = [0,1] \neq \{0\}$, therefore the stronger condition $||\varphi \dotdiv \varphi||_M = ||\varphi||_M \dotdiv ||\varphi||_M$ does not work for $\varphi \to \varphi$, which is a logical tautology, for $||\varphi||_M = [0,1]$. We think that partial models can be helpful in the study of forcing constructions in continuous logic. For model theoretic forcing in continuous logic see [BI].

# Acknowledgements

We are grateful to the anonymous referee for the useful comments and critique of the manuscript, and for bringing [BP] to our attention, in which Ben Yaacov and Pedersen present a direct proof of completeness theorem for **CL** and obtain as a corollary [BP, 9.11] that any computably axiomatizable linear theory is decidable. Unlike [BP], we used an indirect approach of taking a **CL** theory as a theory in **RPL**$\forall$ with some additional axioms and using Hajek's result about completeness of **RPL**$\forall$ to reduce **CL** to **RPL**$\forall$ §2.3 and obtain 3.4. We want to point out that neither proving the completeness of **CL** nor the decidability result was our main intention, rather we used the latter as a stepping stone in proving a constructive version of Henkin's construction in order to build a platform for computable model theory for **CL** similar to the classical case, c.f. [Mil99].

We are in debt to the editors, specially Iraj Kalantari, for their patience and their help in making the paper more readable. Thanks to our colleagues in Domain Theory Group, Amirkabir University of Technology. Some of the results in this paper were presented at IPM Logic Conference 2007, we are thankful to the organizers.